\input amstex
\documentstyle{amsppt} \magnification=1200 \hsize=13.8cm
\catcode`\@=11
\def\NoLogo{\let\logo@\empty}
\catcode`\@=\active
\NoLogo
\def\heat{\lf(\frac{\p}{\p t}-\Delta\ri)}
\def \b {\beta}

\def\lf{\left}
\def\ri{\right}

\def\a{\alpha}

\def\g{\gamma}

\def\p{\partial}

\def\vp{\varphi}

\def\bb{{\bar\beta}}
\def\abb{{\alpha\bar\beta}}
\def\gbd{{\gamma\bar\delta}}

\def \D {\Delta}

\documentstyle{amsppt}
\magnification=1200
\hsize=13.8cm
\vsize=19 cm
\leftheadtext{Huai-Dong Cao and Lei Ni} \rightheadtext{Matrix Li-Yau-Hamilton estimates}
\topmatter
\title{Matrix Li-Yau-Hamilton estimates for the heat equation on K\"ahler manifolds}
\endtitle
\author{Huai-Dong Cao \footnotemark and Lei Ni\footnotemark }
\endauthor
\footnotetext"$^1$"{Research partially supported by NSF grant DMS-0206847.}
\footnotetext"$^2$"{Research partially supported by NSF grant DMS-0203023.}
\address
 {Institute for Pure and Applied Mathematics, Los Angeles, CA 90095-7121 and
Department of Mathematics, Texas A\&M University, College Station, TX 77843}
\endaddress
\email{hcao\@ipam.ucla.edu }
\endemail
\address
{Department of Mathematics, University of California, San Diego, La Jolla, CA 92093}
\endaddress
\email{lni\@math.ucsd.edu}
\endemail
\affil
{Institute for Pure and Applied Mathematics and Texas A\&M University  \\
University of California, San Diego}
\endaffil

\date
August, 2002
\enddate
\endtopmatter
\document

\subheading{ \S1 Introduction} In [L-Y], Peter Li and S.-T. Yau developed the
fundamental gradient estimate, which is now widely called the Li-Yau estimate, for any positive
solution $u(x,t)$ of the heat equation on a Riemannian manifold $M^m$ and showed how the classical
Harnack inequality can be derived from their gradient estimate. When $M^m$ is complete and of
nonnegative Ricci curvature, the Li-Yau estimate is sharp. Later in [H2], Richard Hamilton extended the Li-Yau estimate to the full matrix version of the Hessian estimate of $u$ under the stronger assumptions that $M$ is Ricci parallel and of nonnegative sectional curvature.

In this paper, we consider positive solutions to the heat equation
$$ \heat u(x,t) =0 \tag 1.1 $$
on a complete K\"ahler manifold $M^m$ of complex dimension $m$ with K\"ahler metric
$g=(g_{\abb})$. Our main result is the following complex analogue of Hamilton's Hessian estimate for any positive solution $u$ to (1.1):

\proclaim{Theorem 1.1} Let $M^m$ be a complete K\"ahler manifold of complex dimension $m$ with nonnegative holomorphic bisectional curvature, and $u(x,t)$ be a positive solution of (1.1).
Then, for any vector field $V=(V^\a)$ of type (1,0) on $M$ and $t>0$, we have
$$ u_{\abb}+u_\a V_{\bb}+u_{\bb}V_\a+uV_{\a}V_{\bb}+\frac{u}{t}g_{\abb} \ge 0. \tag 1.2 $$
\endproclaim

A similar result with error terms can be formulated for complete K\"ahler manifolds $M$ with curvature bounded from below.  Note that we have two advantages of being in the K\"ahler category here. Namely, not only we can replace the assumption of nonnegativity of the sectional curvature in Hamilton's result by that of the holomorphic bisectional curvature, but also we can remove the assumption of the Ricci tensor being parallel which is a rather restrictive condition and is definitely needed in Hamilton's work [H2]. Therefore, our Theorem 1.1 should be more applicable.

If we choose the optimal $V=-\nabla u/u$ and take the trace in (1.2), we obtain the gradient estimate of Li-Yau in the K\"ahler case:

$$u_t-\frac{|\nabla u|^2}{u}+\frac{m}{t}u\ge 0. \tag 1.3 $$

We remark that the trace Li-Yau estimate (1.3) is not entirely a corollary of our Li-Yau-Hamilton matrix estimate (1.2) since their result can be obtained under the condition of nonnegative Ricci curvature,
 which is weaker than nonnegativity of the  holomorphic bisectional curvature. However, the conclusion in our result is stronger and therefore the result is more powerful in applications. For example, as pointed out in [H2], the trace Li-Yau estimate (1.3) allows comparisons of the values of $u$ between different points at different times, while the matrix Li-Yau-Hamilton estimate (1.2) also allows the comparisons between different points at the same time.

An immediate application of Theorem 1.1 is the following complex Hessian comparison theorem for the distance function on a complete K\"ahler manifolds of nonnegative holomorphic bisectional curvature:

\proclaim{Corollary 1.1} Let $M$ be a complete K\"ahler manifold with nonnegative holomorphic bisectional curvature. Let $r(x)$ be the distance function to a fixed point $o\in M$. Then in the sense of currents, we have
$$(r^2)_{\abb}(x) \le g_{\abb}(x). \tag 1.4 $$
In particular, when $M$ is noncompact, Buseman functions with respect to geodesics are plurisubharmonic.
\endproclaim

\demo{Proof of Corollary 1.1} Applying Theorem 1.1 to the heat kernel $H(x,y, t)$ with
$V=-\nabla H/H$, we have
$$\left(\log H\right)_{\abb}+\frac{1}{t}g_{\abb}\ge 0.$$
Now it is well known that $H$ is positive and $-t\log H\to r^2(x,y)$ as $t\to 0$. The result
then follows.
\enddemo

Such a Hessian comparison theorem seems to be elusive from the literature even though Greene-Wu (cf [G-W]) have proved the plurisubharmonicity of the Buseman function on such manifolds. We should  mention that Corollary 1.1 was also proved by Jiaping Wang [W] using a different and more direct method earlier.

Now we turn our attention to the trace and matrix estimates for the potential function of the K\"ahler-Ricci flow on a compact or complete noncompact K\"ahler manifold. In the study of the K\"ahler-Ricci flow
$$ \frac{\p}{\p t} g_{\abb}(x,t)=-R_{\abb}(x,t), \tag 1.5$$
it is often useful to consider the time-dependent heat equation:
$$\left(\frac{\p}{\p t}-\Delta_t\right)u(x,t)=0. \tag 1.6 $$
Here $\Delta_t$ denotes the Laplace operator with respect to the evolving metric $g_{\abb}(x,t)$ at time $t$. For example, when $M$ is compact and the first Chern class $c_1(M)=0$, the K\"ahler-Ricci flow was studied by the first author in [C1]. In this case, (1.5) can be reduced to the following scalar complex Monge-Amper\'e flow of the (unknonw) function $\vp(x,t)$:
$$\frac{\p}{\p t}\vp(x,t) = \log\frac{\det(g_{\gbd}(x,t))}{\det(g_{\gbd}(x,0))}+f(x),  \tag 1.7 $$
where $g_{\abb}(x,t)=g_{\abb}(x,0)+\vp_{\abb}(x,t)$ and $f_{\abb}=-R_{\abb}(x,0).$  It is then easy to check that $u=-\vp_t$ is a potential function of the evolving Ricci tensor $R_{\abb}(x,t)$, i.e., the complex Hessian $u_{\abb}(x,t)=R_{\abb}(x,t)$, and satisfies the heat equation (1.6). It is often useful to obtain gradient estimate for positive solutions of (1.6) in general, and in particular for the potential functions of the evolving Ricci tensor. It turns out that the trace estimate of Li-Yau always holds for the positive potential functions of the evolving Ricci tensor without any assumptions on the sign of curvature:

\proclaim{Theorem 1.2} Let $M^m$ be a compact K\"ahler manifold with $c_1(M)=0$. Let $\vp$ and $u$ be given as above, and assume $u>0$. Then we have, for $t>0$,
$$ u_t-\frac{|\nabla u|^2}{u}+\frac{m}{t}u \ge 0.\tag 1.8$$
\endproclaim

A similar result also holds for the $c_1(M)>0$ case. See the statement of Theorem 2.1 in next section.

If we assume that $M$ is complete noncompact with nonnegative holomorphic bisectional curvature, then the matrix gradient estimate holds for the positive potential function $u$. (The compact case would not be of much interests since under the assumptions $c_1(M)=0$ and $M$ has nonnegative holomorphic bisectional curvature, $M$ is in fact holomorphically isometric to a flat complex tori.) Namely, we have

\proclaim{Theorem 1.3} Let $g_{\abb}(t)$, $0\leq t<T$, be a complete solution to the K\"ahler-Ricci flow (1.5) on a noncompact complex manifold $M$ with nonnegative holomorphic bisectional curvature. Let $u(x,t)$ be a positive potential function of the evolving Ricci tensor. Then $u$ satisfies the matrix estimate (1.2).
\endproclaim

In fact, we prove a matrix gradient estimate for any positive solution $u$ to the heat equation
(1.6) coupled with the K\"ahler-Ricci flow (1.5), provided $u$ is plurisubharmonic. See Theorem 3.1 in the last section.

\medskip

\noindent{\bf Acknowledgment.} The authors would like to thank the National Center for Theoretical Sciences at National Tsing Hua University in Hsinchu, Taiwan for the hospitality provided during the writing of this paper. Part of the work was carried out while the first author was visiting the Mathematics Department of Harvard University. He would like to thank Professor S.-T. Yau for making the visit possible and for his encouragement.

\medskip

\subheading{\S2 The Compact Case}
\medskip
Throughout this section, we assume that $M^m$ is compact so that one can apply the tensor maximum principle of Hamilton in [H1] without worrying about any growth assumption on the tensor. We shall first present the proof of Theorem 1.1 in the compact case, and then the proof of Theorem 1.2 as well as the analogous case of $c_1(M)>0$. The proof of Theorem 1.1 in the noncompact case will be given in Section 3.

\demo{Proof of Theorem 1.1 (The compact case)} As in [H2], it suffices to prove that for $t>0$, the Hermitian symmetric (1,1) tensor
$$N_{\abb}:=u_{\abb}+\frac{u}{t}g_{\abb}-\frac{u_{\a}u_{\bb}}{u}\ge 0.\tag 2.1$$
As always, we first apply the heat operator to the tensor $N_{\abb}$. From direct calculations (cf. Lemma 2.1 in [N-T1]), we have
$$\heat u_{\abb}=R_{\abb\gbd}u_{\delta\bar{\gamma}}-\frac{1}{2}R_{\a\bar{s}}u_{s\bb} -\frac{1}{2}R_{s\bb}u_{\a\bar{s}}.\tag 2.2$$
Using the fact that $\D=\frac{1}{2}\left(\nabla_{s}\nabla_{\bar{s}}+\nabla_{\bar{s}}\nabla_s\right)$, we also obtain
$$\split \heat\left(-\frac{u_\a u_{\bb}}{u}\right) & =  \frac{1}{u}\left(
u_{\a s}u_{\bar{s}\bb}
+u_{\a\bar{s}}u_{s\bb}\right) +\frac{2}{u^3}u_{\a}u_{\bb}|u_s|^2 \\
& \ \ +\frac{1}{2u}(R_{\a \bar{s}}u_s u_{\bb}+R_{s\bb}u_\a u_{\bar{s}})\\
& \   \ -\frac{1}{u^2}\left(u_{\a s}u_{\bar{s}}u_{\bb}+u_{\bb s}u_\a u_{\bar{s}}+u_{\a\bar{s}}u_\bb u_{s}+u_{\bb \bar{s}}u_\a u_s\right),
\endsplit \tag 2.3$$
and
$$\heat \left(\frac{u}{t}g_\abb\right)=-\frac{u}{t^2}g_{\abb}.\tag 2.4$$

Combining (2.2)-(2.4), we have
$$\split\heat N_{\abb}& =R_{\abb\gbd}N_{\delta\bar{\g}}
-\frac{1}{2}(R_{\a\bar{s}}N_{s\bb}+N_{\a\bar{s}}R_{s\bb})+
\frac{1}{u}N_{\a\bar{s}}N_{s\bb}-\frac{2}{t}N_{\abb}\\
& \ \ +\frac{1}{u}\left(u_{\a s}- \frac{u_\a u_s}{u} \right)
\left(u_{\bar{s}\bb}-\frac{u_{\bar{s}}u_\bb}{u}\right)+\frac{1}{u}R_{\abb\gbd}
u_{\delta}u_{\bar\g}.
\endsplit \tag 2.5$$

Now according to the tensor maximum principle of Hamilton in [H1], to prove $N_{\abb}\ge 0$ it suffices to show that the right hand side of (2.5) is nonegative when applied to any null vector of $N_{\abb}$. However, it is easy to check that in fact each term on the right hand side of
(2.5) is nonnegative when evaluated at any null vector of $N_{\abb}$. Thus the proof of Theorem 1.1 is proved in the case of $M$ being compact.
\enddemo

\demo{Proof of Theorem 1.2} As in [L-Y], let $v=\log u$. Define $G=t\lf(|\nabla v|^2-v_t\ri)$. It suffices to show that $G\le m$. Direct calculations show that
$$\Delta_t v -v_t =-|\nabla v|^2, \tag 2.6$$
$$\Delta_t|\nabla v|^2=|v_{\a \gamma}|^2 +|v_{\a \bar{\gamma}}|^2
+(\Delta_t v)_{\a}v_{\bar{\a}}+v_{\a}(\Delta_t v)_{\bar{\a}}+R_{\abb}v_{\a}v_{\bar{\b}},
\tag 2.7$$
and $$\frac{\p}{\p t}|\nabla v|^2= R_{\abb}v_{\a}v_{\bar{\beta}}+(v_t)_{\a}v_{\bar{\a}}+ v_{\a}(v_t)_{\bar{\a}}. \tag 2.8$$ Here the first term on the right hand side of (2.8) is due to the fact that we have a time-dependent metric evolving by the K\"ahler-Ricci flow (1.5).

From (2.6) we also have
$$v_{tt}-\Delta_t(v_t)=R_{\abb}(v_{\beta\bar{\a}}+v_{\beta}v_{\bar{\a}})
+(v_t)_{\a}v_{\bar{\a}}+v_{\a}(v_t)_{\bar{\a}}. \tag 2.9$$

From (2.7)--(2.9) it follows
$$\split\left(\Delta_t -\frac{\p}{\p t}\right)(|\nabla v|^2-v_t) & = |v_{\a \gamma}|^2 +|v_{\a \bar{\gamma}}|^2 +R_{\abb}(v_{\beta\bar{\a}}+v_{\beta}v_{\bar{\a}})\\
&\ \ \ -(|\nabla v|^2-v_t)_{\a}v_{\bar{\a}}-(|\nabla v|^2-v_t)_{\bar{\a}}v_{\a}\\
& \ge |v_{\a \bar{\gamma}}|^2 -(|\nabla v|^2-v_t)_{\a}v_{\bar{\a}} -(|\nabla v|^2-v_t)_{\bar{\a}}v_{\a}.
\endsplit
\tag 2.10$$ Here we have used the fact that $R_{\abb}(v_{\beta\bar{\a}}+v_{\beta}v_{\bar{\a}})=R_{\abb}u_{\beta\bar{\a}}/u=|R_{\abb}|^2/u \ge 0$.

From (2.10), we obtain
$$\split\left(\Delta_t-\frac{\p}{\p t}\right)G &\ge t|v_{\a \bar{\gamma}}|^2 -2<\nabla G, \nabla v> -\frac{G}{t}\\
& \ge \frac{t}{m}\left(\Delta_t v\right)^2 -2<\nabla G, \nabla v>-\frac{G}{t}\\
& = \frac{G^2}{tm}-2<\nabla G, \nabla v>-\frac{G}{t}.
\endsplit
\tag 2.11$$

Applying the maximum principle argument to the above inequality, it then follows that $G\le m$, which completes the proof of the theorem.
\enddemo

In the case of compact $M$ with first Chern class $c_1(M)>0$, we can obtain a similar result. In this case, consider the normalized K\"ahler-Ricci flow
$$ \frac{\p}{\p t} g_{\abb}(x,t)=-R_{\abb}(x,t) +g_{\abb}(x,t) \tag 2.12$$
with the initial metric $g(x,0)$ and its K\"ahler form $\omega$ such that $c_1(M)=\pi [\omega]$. Similar to (1.5), (2.12) can also be reduced to a complex Monge-Amper\'e flow of the form $$\frac{\p}{\p t}\vp(x,t) = \log\frac{\det(g_{\gbd}(x,t))}{\det(g_{\gbd}(x,0))}+\vp(x,t)+f(x).  \tag 2.13 $$
Here again $g_{\abb}(x,t)=g_{\abb}(x,0)+\vp_{\abb}(x,t)$ and $f_{\abb}=g_{\abb}(x,0)-R_{\abb}(x,0).$  Furthermore, it was shown by the first author [C1, C2] that the solution to (2.13), hence also (2.12), exists for all time.

Set $w=-\vp_t$. Then $w$ is a potential function satisfying
$$w_{\abb}(x,t)=R_{\abb}(x,t)-g_{\abb}(x,t) \tag 2.14$$
and
$$\lf(\frac{\p}{\p t}-\Delta_t\ri) w =w. \tag 2.15$$

Similar to Theorem 1.2, we have
\proclaim{Theorem 2.1} Let $M^m$ be a compact K\"ahler manifold with $c_1(M)>0$. Let $\vp$ and $w$ be defined as above, and assume $w>0$. Then we have, for $t>0$,
$$ w_t-\frac{|\nabla w|^2}{w}+\frac{m}{t}w \ge w>0.\tag 2.16$$
\endproclaim

\demo{Proof of Theorem 2.1} Let $u=e^{-t}w$, then $u$ is a positive solution to the heat equation (1.6) coupled with (2.12). As in the proof of Theorem 1.2, we let $v=\log u$ and define $G=t\lf(|\nabla v|^2-v_t\ri)$. It follows from similar calculations there that
$$\split\left(\Delta_t -\frac{\p}{\p t}\right)(|\nabla v|^2-v_t) & = |v_{\a \gamma}|^2 + |v_{\a \bar{\gamma}}|^2 +\left(R_{\abb}-g_{\abb}\right)(v_{\beta\bar{\a}}+v_{\beta}v_{\bar{\a}})\\
&\ \ \ -(|\nabla v|^2-v_t)_{\a}v_{\bar{\a}}-(|\nabla v|^2-v_t)_{\bar{\a}}v_{\a}+
|\nabla v|^2\\
& \ge  |v_{\a \bar{\gamma}}|^2-(|\nabla v|^2-v_t)_{\a}v_{\bar{\a}} -(|\nabla v|^2-v_t)_{\bar{\a}}v_{\a}.
\endsplit
\tag 2.17
$$
Here we have used the fact that
$$(R_{\abb}-g_{\abb})(v_{\beta\bar{\a}}+v_{\beta}v_{\bar{\a}})=(R_{\abb}-g_{\abb})u_{\beta\bar{\a}}/u=|R_{\abb}-g_{\abb}|^2/w \ge 0.$$
Hence $G$ satisfies the same differential inequality (2.11), and we can conclude the same way
that $G\le m$. Therefore, the function $u=e^{-t}w$ satisfies the estimate (1.8). Expressing this in terms of $w$, we obtain the desired estimate (2.16) and thus proves Theorem 2.1.
\enddemo



\medskip
\subheading{\S 3 The Complete Noncompact Case}
\medskip
Now we consider the case when $M$ is a complete noncompact K\"ahler manifold with nonnegative holomorphic bisectional curvature. Due to the fact that uniqueness of the solution to the scalar heat equation fails to be true in general on a noncompact manifold, one normally needs to impose some kind of growth conditions on the function $u$ as well as its first and second order derivatives in order to be able to apply Hamilton's tensor maximum principle (or its argument) to the tensor $N_{\abb}$ defined in (2.1). However, in our case of proving Theorem 1.1, we shall see that we can get away without imposing any growth assumptions on $u$ and its derivatives. The key here is that we are working with a positive solution of the heat equation, thus we can make use of the available estimate of Li-Yau to obtain the required growth estimates at any positive time. First let us collect some basic facts.

\proclaim{Lemma 3.1} Let $u(x,t)$ be a positive solution to (1.1). Then we have
$$\heat |\nabla u|^2 \le -\|u_{\abb}\|^2-\|u_{\a\beta}\|^2 \tag 3.1$$
and
$$\heat \|u_{\abb}\|^2 \le 0. \tag 3.2$$
\endproclaim

\demo{Proof} Both (3.1) and (3.2) can be verified by direct calculaitons. Here
the nonnegativity of the holomorphic bisectional curvature is needed. For more details, see for example Lemma 1.1 in [N-T1] and Lemma 1.5 in [N-T2].
\enddemo

We also need to use the result of Li-Yau on the Harnack inequality for positive solution to the heat equation. Let $o\in M$ be a fixed point, and let $u(x,t)$ be a positive solution of (1.1).
Since our focus here is to obtain a upper bound on $u$ for positive time we can assume, without the lose of generality, that $u(x,t)$ is defined on $M\times [0,2]$. By the Harnack inequality in [L-Y] (Theorem 2.2(i), page 168) we have, for $0<t<1$
$$u(x,t)\le \frac{C}{t^m}u(o,2)\exp(ar^2(x)). \tag 3.3$$
Here $a>0$ is a constant and $r(x)$ is the distance function from the point $o$. In particular, for $t\ge \delta>0$, there exists a constant $b>0$  (might depends on $\delta$) such that
$$u(x,t)\le \exp(b(r^2(x)+1)). \tag 3.4$$
In fact using (3.1) and (3.2) together with the mean value inequality of Li-Tam  we can push further to obtain the similar control on the gradient and the complex Hessian of $u$ for $t>2\delta$. In fact, we have the following

\proclaim{Lemma 3.2} For $t\ge 2\delta$, there exists $b_1\ge b>0$ such that
$$|\nabla u|^2(x,t)\le \exp(b_1(r^2(x)+1)) \tag 3.5$$
and $$\|u_{\abb}\|^2(x,t)\le \exp(b_1(r^2(x)+1)). \tag 3.6 $$
\endproclaim

\demo{Proof} First we prove that for some $b_2>0$,
$$\int_{\delta}^{T}\int_M \exp(-b_2(r^2(x)+1))|\nabla u|^2(x)
\, dx\, dt < \infty \tag 3.7$$
To see this, we multiply $\vp^2$, where $\vp$ is a cut-off function, on both sides of the equation
$$\heat u^2 =-2|\nabla u|^2$$
and then integrate by parts. As in [N-T2] we have
$$\split 2\int_0^T\int_M \vp^2|\nabla u|^2 dxdt&=-\int_0^T\int_M \vp^2\heat u^2\\
&\le \int_M \vp^2 u_0^2(x) dx+4\int_0^T\int_M  \vp u|\nabla \vp|\,|\nabla u|dxdt\\
&\le \int_M \vp^2 u_0^2(x) dx+4\int_0^T\int_M |\nabla \vp|^2u^2dxdt+\int_0^T\int_M \vp^2|\nabla u|^2 dxdt.
\endsplit$$

Now (3.7) follows from (3.3). As in the proof of Lemma 1.6 in [N-T2], applying the mean value inequality of Li-Tam (Theorem 1.1 of [L-T]) and using (3.1) imply (3.5). Similar argument using (3.1) and (3.2) proves (3.6). For more details, see the proof of Lemma 1.6 in [N-T2].
\enddemo

Now we are in the position to prove Theorem 1.1 for the complete noncomapct case.

\demo{Proof of Theorem 1.1 (The noncompact case)} We first shift the time by $2\delta$. By doing so, $u(x,t)$ togeter with its gradient and complex Hessian satisfy (3.3)--(3.5). If we can prove the theorem for this case, then we would have (1.2) when replacing $t$ by $t+2\delta$. By letting $\delta\to 0$ we would complete the proof of Theorem 1.1. Therefore, without loss of generality, we can assume (3.3)--(3.5) hold. By a similar argument we also can assume $u\ge \delta$ in the proof.

As in the proof of Lemma 2.1 in [N-T2], we first construct a function $\phi(x,t)$ such that
$$\heat \phi =\phi$$
and
$$\phi(x,t)\ge C_1\exp(2b_1(r^2(x)+1))$$
for some constant $C_1>0$. Let $N_{\abb}$ be the Hermitian (1,1) tensor defined in (2.1). We consider the (1,1) tensor $Z_{\abb}=t^2N_{\abb}+\epsilon \phi g_{\abb}$, where $g_{\abb}$ is the metric tensor. Clearly we only need to show that $Z_{\abb}\ge 0$ for any $\epsilon >0$. We shall prove this by contradiction. Suppose it is not true, then by the growth nature of $\phi$ and the fact that $N_{\abb}>0$ at time $t=0$,  we know that there exists a first time $t_0>0$, and a  point $x_0\in M$
and a unit vector $V=v^{\a}\frac{\p}{\p z_{\a}}\in T_{x_0}M$  such that
$Z_{\abb}(x_0, t_0)v^\a\bar{v}^{\beta}=0$.
Now we choose a normal coordinate around $x_0$ and extend $V$ to be a local unit vector field near $x_0$ by parallel translation along the geodesics emanating from $x_0$. It then follows from the direct calculation that, at point $x_0$,
$$\Delta \lf(Z_{\abb}v^\a\bar{v}^\beta\ri)=\lf(\Delta Z_{\abb}\ri)v^{\a}\bar{v}^\beta.$$
Since $Z_{\abb}v^\a\bar{v}^\beta\ge 0$ for all $(x,t)$ with $t\le t_0$ and $x$ close to $x_0$, and $Z_{\abb}v^\a\bar{v}^\beta=0$ at $(x_0, t_0)$ we see that at $(x_0, t_0)$,
$$0\ge \heat \lf(Z_{\abb}v^\a\bar{v}^\beta\ri).$$
On the other hand, using (2.5) we also have, at $(x_0, t_0)$,
$$\split\heat \lf(Z_{\abb}v^\a\bar{v}^\beta\ri) &= \lf(\heat Z_{\abb}\ri)v^\a\bar{v}^\beta\\
&\ge t^2\lf(R_{\abb\gbd}N_{\bar{\g}\delta}
-\frac{1}{2}R_{\a\bar{s}}N_{s\bb}-\frac{1}{2}R_{s\bb}N_{\a\bar{s}} \ri) v^\a\bar{v}^\beta\\
& \ \ \ + \frac{t^2}{u} R_{\abb\gbd}u_{\bar{\g}}u_{\delta} v^\a\bar{v}^\beta+
\epsilon \phi |V|^2\\
& \ge R_{\abb\gbd}Z_{\bar{\g}\delta}v^\a\bar{v}^\beta -\frac{1}{2}R_{\a\bar{s}}Z_{s\bar{\beta}}
v^\a\bar{v}^\beta -\frac{1}{2}R_{s\bb}Z_{\a\bar{s}}v^\a\bar{v}^\beta +\epsilon \phi |V|^2\\
&>0.
\endsplit$$

We now have a  contradiction. Therefore we have completed the proof Theorem 1.1 in case $M$ is noncompact.
\enddemo

\proclaim{Remark} As we mentioned before, our result does not contain Li-Yau's gradient estimate even though by taking trace we obtain their estimate (1.3) since we need to use the Harnack inequality and the mean value inequality of Li-Tam in Lemma 3.2. The proof of both the Harnack inequality and Li-Tam's mean value inequality rely on Li-Yau's gradient estimate. Also we need to assume nonnegativity of the holomorphic bisectional curvature in stead of the Ricci curvature. On the other hand our result is about the full Hessian matrix of the function and therefore is stronger than Li-Yau's trace estimate.
\endproclaim

Finally, we consider the matrix gradient estimate for any positive solution $u$ to the heat equation (1.6) coupled with the K\"ahler-Ricci flow (1.5). In this case, (2.2) and (2.3) remain the same but (2.4) and (2.5) become, respectively,
$$\heat\left(\frac{u}{t}g_\abb\right)=-\frac{u}{t^2}g_{\abb}-\frac{u}{t}R_{\abb}.\tag 2.4\'$$
and
$$\split\heat N_{\abb}& =-\frac{1}{2}R_{\a\bar{s}}N_{s\bb}-\frac{1}{2}N_{\a\bar{s}}R_{s\bb}+
\frac{1}{u}N_{\a\bar{s}}N_{s\bb}-\frac{2}{t}N_{\abb}\\
& \ \ +\frac{1}{u}\left(  u_{\a s}- \frac{u_\a u_s}{u} \right)
\left(u_{\bar{s}\bb}-\frac{u_{\bar{s}}u_\bb}{u}\right)+R_{\abb\gbd}
u_{\bar{\g}\delta}.\endsplit \tag 2.5\'
$$
Notice that under the extra assumption that $u(x,t)$ is plurisubharmonic, the last term is nonnegative definite. Therefore we can prove a similar result as in Theorem 1.1 for the coupled case:

\proclaim{Theorem 3.1} Let $g_{\abb}(x,t)$ be complete K\"ahler metrics evolving by the K\"ahler-Ricci flow (1.5) on $M^m$, and $u(x,t)$ be a positive solution to the time-dependent heat equation (1.6). Assume that the holomorphic bisectional curvature of $g_{\abb}(x,t)$ is nonnegative and $u(x,t)$ is plurisubharmonic. Then
$$u_{\abb}+u_\a V_{\bb}+u_{\bb}V_\a+uV_{\a}V_{\bb}+\frac{u}{t}g_{\abb}\ge 0.$$
\endproclaim

\demo{Proof of Theorem 1.3} Apply Theorem 3.1 to the potential function $u$, which is plurisubharmonic since its complex Hessian is equal to the Ricci tensor.
\enddemo

Notice in [N-T3], the authors proved that under some average curvature decay assumption one indeed can obtain the potential function $u(x,t)$ for the Ricci tensor by solving the Poincar\'e-Lelong equation and utilizing the volume element.

\proclaim{Remarks} (i) Taking the trace in Theorem 3.1, we can obtain the gradient estimate for $u$ obtained before in [N-T1]. Notice again that this gradient estimate in [N-T1] is not entirely a corollary of Theorem 3.1 since there they only need to assume that the Ricci curvature is nonnegative while in Theorem 3.1 we need to assume that the holomorphic bisectional curvature is nonnegative.

(ii) In [N-T1, N-T2], the authors also studied the question under what conditions the plurisubharmonicity of $u(x,t)$ will be preserved by the heat flow in the time-dependent case.
\endproclaim

\enddocument